\documentclass[10pt,oneside]{article}
\usepackage{tikz}
\usepackage{tikz-cd}
\usepackage{amssymb}
\usepackage{amsfonts}
\usepackage{amscd}
\usepackage{verbatim}
\usepackage{latexsym} 
\usepackage{graphicx}
\usepackage{amsthm}
\usepackage{latexsym}
\usepackage{amsmath}
\usepackage{lscape}
\usepackage{palatino, url, multicol}
\usepackage{subfig}
\usepackage{geometry}               
\theoremstyle{definition}
\newtheorem{theorem}{Theorem}[section]
\newtheorem{lemma}[theorem]{Lemma}
\newtheorem{proposition}[theorem]{Proposition}
\newtheorem{corollary}[theorem]{Corollary}

\theoremstyle{definition}
\newtheorem{remark}[theorem]{Remark}
\theoremstyle{definition}

\theoremstyle{definition}

\theoremstyle{definition}

\theoremstyle{definition}

\title{Devil's Staircase -- Rotation Number of \\ Outer Billiard with Polygonal Invariant Curves }
\author{Zijian Yao }

\begin{document}
\maketitle
\begin{abstract}
\noindent
In this paper, we discuss rotation number on the invariant curve of a one parameter family of outer billiard tables. Given a convex polygon $\eta$, we can construct an outer billiard table $\mathcal{T}$ by cutting out a fixed area $\mathcal{A}$ from the interior of $\eta$. $\mathcal{T}$ is piece-wise hyperbolic and the polygon $\eta$ is an invariant curve of $\mathcal{T}$ under the billiard map $\phi$. We will show that, if $\beta \in \eta$ is a periodic point under  $\phi$ with rational rotation number $\tau = \frac{p}{q}$, then $\phi^{q}$ is not the local identity at $\beta$. This proves that the rotation number $\tau$ as a function of the parameter $\mathcal{A}$ is a devil's staircase function.

\end{abstract}

\section{ \fontsize{15}{14}\selectfont \textit{Introduction}}

\noindent
In mid 1990s, Gutkin and Knill \cite{GE} considered a one parameter family of inner billiard tables which have an equilateral triangle as a common caustic (The billiard tables can be constructed geometrically by the string construction, where the length $l$ of the string is the parameter).  The family of circle homeomorphisms obtained by restricting the billiard map to the canonical invariant circles gives a family of associated rotation numbers with parameter $l$. They proved that the rotation number $\tau(l)$, as a function of $l$, is a devil's staircase function. This means that, in this one parameter family of tables, there does not exist one consisting solely of periodic points. For a concise introductory treatment of rotation numbers, I refer the readers to \cite{KA}. 
However, this phenomenon is not universal. In 1988, Innami \cite{IN} already gave descriptions of a family of smooth inner billiard tables that consist only of 3-periodic points.  
In 2006,  Baryshnikov and Zharnitsky \cite{BY} also studied inner billiard with full one parameter family of periodic orbits. They showed that there exist billiard tables which consist only of periodic points but have no elliptic boundaries. \\

\noindent
In this article we study a related problem on outer billiard systems. Introductory treatments on outer billiards can be found in \cite{DT, GD, TS}. \\

\noindent
Let $D$ be an outer billiard table and $C$ an invariant curve of the outer billiard map $\phi$. Let $x$ be a point on $C$ and $y $ its image $\phi(x)$. It is known that the area bounded by the segment $xy$ and the invariant curve $C$ is constant for all $x \in C$. In other words, $D$ can be recovered from $C$ as an envelope of segments of constant area ( see \cite{TS} for more details ).  Now we construct the outer billiard tables by cutting a fixed area from a convex polygon. It is known that the area construction results in a piecewise hyperbolic table, which has the original polygon as its invariant curve. \\

\noindent
For a simple example we start with the square $P_1 P_2 P_3 P_4$ as shown in Figure \ref{figure01}. We label the vertices and sides as the figure suggests. The resulted table is a symmetric piecewise hyperbolic ``square''. In 2006, Genin \cite{GD2} studied precisely this one parameter family of tables, and showed that the orbits inside the square invariant curve have chaotic behaviours.\\

\begin{center}
\begin{tikzpicture}
[scale=3, vertices/.style={draw, fill=black, circle, inner sep=1pt}]
\node[vertices] (a) at (0,0) {};
\node[vertices] (b) at (1,0) {};
\node[vertices] (c) at (1,1) {};
\node[vertices] (d) at (0,1) {};
\filldraw[fill=white, draw=black] (0,0)--(1,0)--(1,1)--(0,1)--cycle;
\put(-15,-10){$P_1$}
\put(85,-10){$P_2$}
\put(85,90){$P_3$}
\put(-15,90){$P_4$}
\put(35,-10){$side \, 1$}
\put(85,40){$side \, 2$}
\put(35,90){$side \, 3$}
\put(-30,40){$side \, 4$}
\end{tikzpicture} \qquad  \qquad \qquad \qquad \qquad  \begin{tikzpicture}
[scale=3, vertices/.style={draw, fill=black, circle, inner sep=1pt}]
\node[vertices] (a) at (0,0) {};
\node[vertices] (b) at (1,0) {};
\node[vertices] (c) at (1,1) {};
\node[vertices] (d) at (0,1) {};
\node[vertices] (e) at (0.6,0) {};
\filldraw[fill=black!10, draw=black] (0,0)--(0.6,0)--(0,1)--cycle;
\draw (e)--(b);
\draw (c)--(b);
\draw (d)--(c);
\put(20,-10){$a $}
\put(0,25){$A = \frac{a}{2} $} 
\end{tikzpicture} \\
\vspace*{0.1in}
\end{center} 
\captionof{figure}{}
\label{figure01}
\vspace*{0.15in}

Recall that given a circular homeomorphism $f: S^1  \rightarrow S^1$, the natural projection $\pi : \mathbb{R} \rightarrow S^1$ provides a lift of the map $f$ to homeomorphism $F: \mathbb{R} \rightarrow \mathbb{R}$ such that $\pi \circ F = f \circ \pi$. It is known that $F$  is unique up to adding integer constants. The rotation number $\tau$ of the map $f$ is defined as: 
$$\displaystyle \tau _f = \pi \large( \lim_{n \rightarrow \infty} \frac{F^n (x) - x}{n}  \large).$$

The following facts about rotation number are due to Poincar\'{e}. \cite{KA}

\textbf{[Fact 1]} Let $f: S^1 \rightarrow S^1$ and $F: \mathbb{R} \rightarrow \mathbb{R}$ be as above , then the limit defined above exists for all $ x \in \mathbb{R}$. 

\textbf{[Fact 2]} Let $f': X \rightarrow X$ be a homeomorphism, where $X$ is homeomorphic to a circle by $h : X \rightarrow S^1$, then $\tau(f') := \tau (h^{-1} \circ f \circ h) = \tau (f)$.  In particular, the rotation number is independent of the choice of the starting point.

\textbf{[Fact 3]}{\label{001}} $\tau _f = \frac{p}{q} \in \mathbb{Q}$ if and only if $f$ has a periodic orbit of period $q$ (assume $(p, q) = 1$).

\textbf{[Fact 4]} $\tau(\cdot)$ is continuous in the $C^0$ topology. \\

A continuous and non-decreasing function $\varphi: [0, 1] \rightarrow \mathbb{R}$ is called a \textbf{devil's staircase} if there is a family of disjoint open subintervals of $\mathcal{I} = [0,1]$ such that the union of all these subintervals is dense on $\mathcal{I}$ and the function $\varphi$ takes distinct constant value at each of the subintervals. \\ 

The area construction from an arbitrary convex polygon $\eta$ gives a self-map $\phi_a: \eta \rightarrow \eta$, which in turn yields a circular homeomorphism $f_a: S^1 \rightarrow S^1$. In the rest of the article we will not distinguish $f_a$ from $\phi_a$.  Denote its associated rotation number by $\tau (a)$,  as a function of area ($a = 2A$, where $A$ is the area cut off). Our main result is that $\tau(a)$ as a function of $a$ is a devil's staircase function. \\

It is known that the rotation number $\tau$ is increasing at points $a$ when $\tau(a) $ is irrational and is constant at points $a$ if $\tau(a)$ is rational, as long as not all points are periodic under outer billiard map $\phi_a$ (in which case, the circular map is conjugate to a rotation) \cite{KA}. Therefore, it suffices to prove that no $n^{th}$ iteration  $(\phi_a)^n$ is identity on the polygonal invariant curve.

\section{ \fontsize{15}{14}\selectfont \textit{Main Theorem}}

Consider a generic convex polygon $\eta = P_1 P_2 ... P_n$. Vertices and sides of the polygon are labelled as in Figure \ref{figure02}. We follow the area construction to obtain a table for which $\eta$ is an invariant curve. Let $\tau(A) = \frac{p}{q} \in \mathbb{Q}$ and let $\mathcal{O} = \beta_0 \beta_1 ... \beta_q$ be a corresponding $q$-periodic orbit, where $\beta_q = \beta_0$. We know that $\tau^{-1} (\frac{p}{q})$ is a nonempty open interval if not all points are periodic on $\eta$. \\

Our goal is to prove that no $q^{th}$ iteration $\phi_a^q$ can be identity for $\eta$ if $ 0 < A = \frac{a}{2} < \frac{1}{2}$, i.e.,  we aim to prove the following theorem:  \\

\begin{center}
\includegraphics[scale=0.3]{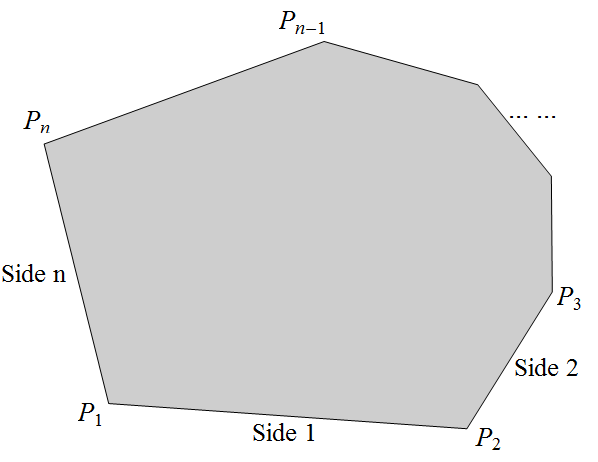} \\
\captionof{figure}{}
\label{figure02}
\end{center} 
\vspace*{0.1in}

\theorem[Main Theorem]{\label{005}}  \par \indent 

Let $\tau(a_0) = \frac{p}{q} \in \mathbb{Q}$ for a convex polygon $\eta$ and $\mathcal{O} = \beta_0 \beta_1 \beta_2 ... \beta_q$ be a corresponding $q$-periodic orbit. Let $\phi_a$ be map $\eta \rightarrow \eta$ given by area construction with area parameter $a = 2A$. Then $\phi_a \, ^q $ is not identity on any non-empty open interval containing $\beta_0$.\\

We shall prove the theorem in the remaining part of this section. First we assume $\mathcal{O}$ does not contain any vertex of $\eta$. We will need to be slightly more careful if there are some $P_j \in \mathcal{O}$ but most arguments still apply.\\

Consider the area cutting line $ l_i  $ defined as the line containing segment $ \beta_i \beta_{i+1}$. Since $\eta$ is convex, $l_i$ intersects $\eta$ at exactly two points, $\beta_i$ and $\beta_{i+1}$.   Let the sides containing $\beta_i$ and $\beta_{i+1}$ be $l_{i1}$ and $l_{i2}$ respectively. Further assume that $l_{i1}$ and $l_{i2}$ are not parallel, so they intersect at some point $Q_{i} $, which clearly does not lie on $l_i$.  $l_i$ divides the plane  $\mathbb{R}^2$ into two open half planes. Let $\mathcal{A}_i$ be the part of interior of $\eta$ that is cut off from the area construction. It is important that we always cut area less than half of the area enclosed by $\eta$, so $S(\mathcal{A}_i) = A < \frac{1}{2} S_{\text{total}} $. \\

\begin{center}
\includegraphics[scale=0.35]{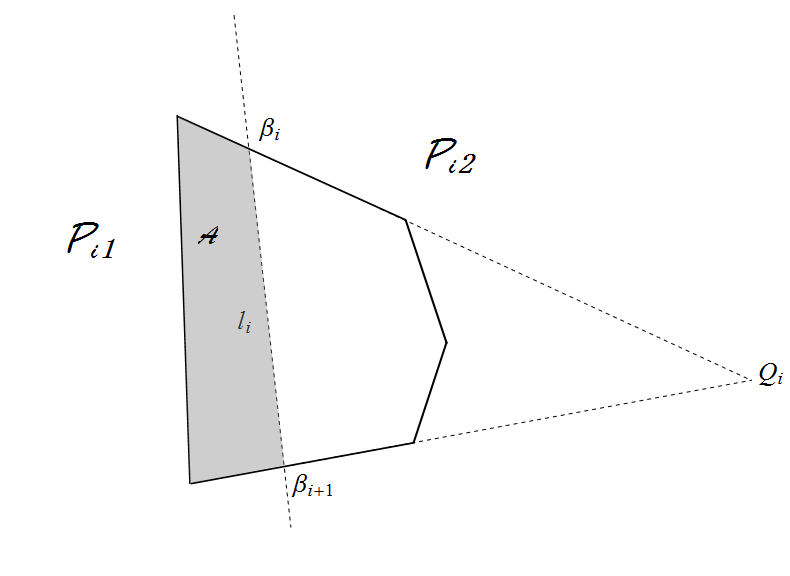} \\
\captionof{figure}{}
\label{figure03}
\end{center}

Define $\mathcal{P}_{i1}$ to be the open half plane containing $\mathcal{A}_i$, and $\mathcal{P}_{i2} = \mathbb{R}^2 - l_i - \mathcal{P}_1 $ to be the other open half plane : [Figure \ref{figure03}]. \\

Note that the definition of $\mathcal{P}_{i1}$ and $\mathcal{P}_{i2}$ varies for different $l_i$, so the sub-index $i$ is necessary to distinguish the division for each line $l_i$. 
Since $Q_i \notin l_i$, then either $Q_i \in \mathcal{P}_{i1}$ or  $Q_i \in \mathcal{P}_{i2}$ . We say that the line $l_{i}$ is \textbf{good} if $Q_i \in \mathcal{P}_{i1}$.
Set $\sigma (l_i) = \sigma (\beta_i \beta_{i+1}) = 1$ if $l_i$ is good. On the other hand, if $Q_i \notin \mathcal{P}_{i1}$, set $\sigma (l_i) =\sigma (\beta_i \beta_{i+1})= -1$.  We say $\sigma (l_i) =\sigma (\beta_i \beta_{i+1})= 0$ if the two sides that $l_i$ intersects are parallel. Figure \ref{figure04} illustrates the definitions. \\

\begin{center}
\includegraphics[scale=0.35]{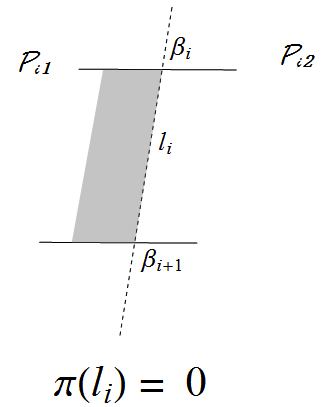} 
\includegraphics[scale=0.35]{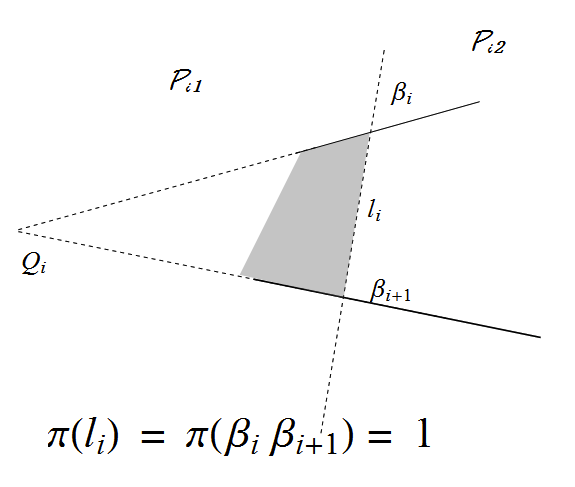}  \\
\includegraphics[scale=0.35]{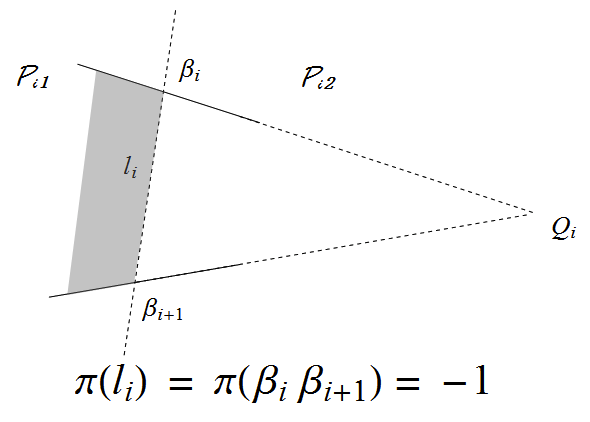} \\
\captionof{figure}{}
\label{figure04}\end{center} 
\vspace*{0.2in}

\lemma {\label{001}} \indent 

If  $\sigma (l_i) = \sigma (\beta_i \beta_{i+1})= -1$, then $\sigma (l_{i+1}) = \sigma ({\beta_{i+1} \beta_{i+2}})= 1$  and $\sigma (l_{i-1}) = \sigma ({\beta_{i-1} \beta_{i}})= 1$. 
\begin{proof}

$l_i$ intersects with two sides $\beta_i Q_i$ and $\beta_{i+1} Q_i$, where $Q_i$ is the intersection of the two sides as defined above. As Figure  \ref{figure05} illustrates, $Q_i \in \mathcal{P}_{i2}$.  Since the area cut off is strictly less than a half, we know that $\beta_{i+2} \in \mathcal{P}_{i2}$, and the area $\mathcal{A}_{i+1}$ cut off from the line $l_{i+1}$ is also in part $\mathcal{P}_{i2}$. Thus $\mathcal{A} _{i+1} \subset \mathcal{P}_{i2}$. The new intersection $Q_{i+1}$, if exists, necessarily lies on $\beta_{i+1} Q_i$. Thus we know $\sigma (l_{i+1}) = 1$ if and only if $Q_{i+1}$ lies on the right hand side of $\beta_{i+1}$, i.e., $Q_{i+1} \in \mathcal{P} _{i2}$; \quad $\sigma (l_{i+1})= -1$ if and only if $Q_{i+1} \in \mathcal{P} _{i1}$; \quad and if $Q_{i+1}$ does not exist, then $\sigma (l_{i+1}) = 0$. \\

Assume for contradiction that $\sigma (l_{i+1}) \neq 1$. First assume $\sigma (l_{i+1}) = -1$, so $Q_{i+1} \in \mathcal{P} _{i1}$. Since $\eta$ is convex, by elementary geometry we know that $\eta$ is contained in the triangular wedge $\angle Q_{i+1} Q_i \beta_i$ and $\angle Q_i Q_{i+1} \beta_{i+2}$. So $\eta \in \triangle Q_{i+1} Q_i Q'$ where $Q'$ in the intersection of $Q_i \beta_i$ and $Q_{i+1} \beta_{i+2}$, as illustrated in Figure  \ref{figure05}.

\begin{center}
\includegraphics[scale=0.4]{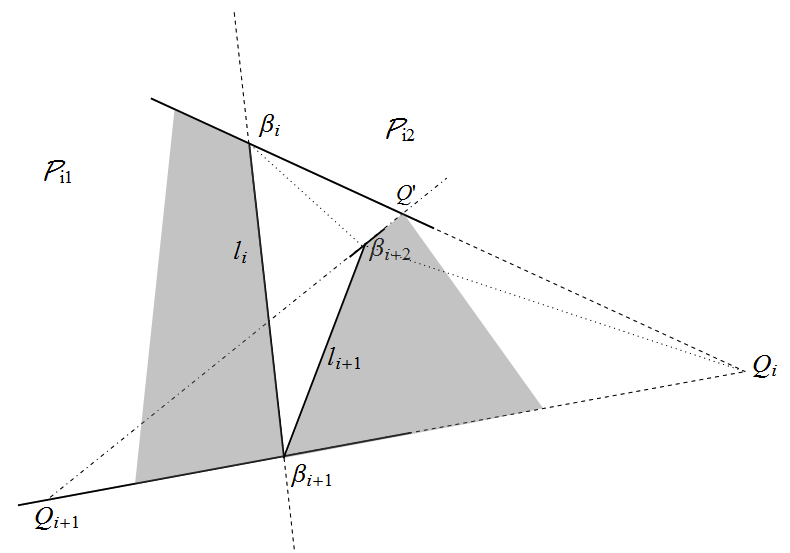}\\
\captionof{figure}{}
\label{figure05}
\end{center} 
\vspace*{0.1in}

This is a contradiction since $\eta$ contains points that lie on the left hand side of $\beta_i$ on the extension of $Q_i\beta_i$. So $\sigma (l_{i+1}) \neq -1$. 

If $\sigma(l_{i+1}) = 0$, then we get a similar contradiction based on geometric argument. \\

The proof of the second part of the lemma uses a similar argument, going in the opposite direction from $l_i$. 
\end{proof}

Back to our $q$-periodic orbit $\mathcal{O} = \beta_0 \beta_1 ... \beta_q$ on $\eta$. As illustrated in Figure  \ref{figure06}, we define a sequence of lengths $a_1,   b_1, a_2, b_2,  ...,   a_l, b_l$ as the following:  we start from the point $\beta_0$, find the smallest index $i_1, \, \, 0 \leq i_1 \leq q , \, \, $ such that $\beta_{i_1}$ and $\beta_{i_1+1}$ are not on parallel sides. Define $a_1$ to be the length from $\beta_{i_1}$ to the vertex to its right (here by right I mean the adjacent vertex in counterclockwise orientation), i.e., if $\beta_{i_1}$ is on side $s, \, \, s \in \{1, 2, ..., n \}$, then $a_1$ is the length between $\beta_i$ and $P_{s+1} $ (whenever the subscripts exceed $n$, we reduce $\mod n$). Define $b_1$ to be the length between $\beta_{i_1 + 1}$ and the vertex to its left, so if $\beta_{i_1 + 1}$ is on side $s'$, then $b_1 = \beta_{i_1+1} P_{s'}$. Now if we extend sides $s$ and $s'$, they necessarily intersect at some point $Q_{i_1}$.  Let the length $Q_{i_1} P_{s+1}$ be $d_1$ and the length $Q_{i_1} P_{s'}$ be $d_1'$. Note that $P_{s+1}, P_{s'}, Q_i$ could coincide, in which case $d_1 = d_1' = 0$.  \\

Next, we start from the point $\beta_{i_1+1}$ and repeat the process above to define $a_2, \, b_2, \, ... , a_l, b_l $. We also define the intersections $Q_{i_1}, Q_{i_2}, ..., Q_{i_l}$ and lengths $d_1, d_1', d_2, d_2', ... ,d_l, d_l'$ analogously.

\lemma {\label{002}}\par \indent

(1). If $\sigma (\beta_{i_k} \beta_{i_k+1}) = 1$, then $b_k + d_k' = \displaystyle {{c_0}}/{(a_k + d_k)}$ for some constant $c_0$ depending on $A$ and $\eta$;

\quad \,\,\, If $\sigma (\beta_{i_k} \beta_{i_k+1}) = -1$, then $ d_k' - b_k = \displaystyle  {{c_0'}}/({d_k - a_k})$ for some other constant $c_0'$. 

(2).  There exists a constant $c_k'$ such that $a_{k+1} + b_k = c_k'$ regardless of whether $l_{i_k}$ is good or not. \\

\begin{proof}
We prove the special case when $k =1$. Denote the area of a polygon $\xi$ by $S(\xi)$.  As Figure  \ref{figure06} indicate below, if $l_{i_1} = \beta_{i_1} \beta_{i_1+1}$ is good, then $$ \frac{1}{2} (b_1 + d_1')(a_1 + d_1) \sin \angle Q_{i_1}  =S (\triangle Q_{i_1} \beta_{i_1}\beta_{i_1+1} ) = A + S (Q_{i_1} P_{s+1} P_{s+2} ...  P_{s'}  ) = \text{constant}.$$ Thus the product $(b_1 +d_1') (a_1 + d_1)$ is a constant. The proof of the case when $l_{i_1}$ is bad is similar, as shown in Figure \ref{figure06}.  

We know that if $\sigma (\beta_{i_1} \beta_{i_1+1}) = \pm 1$, then $b_1 + a_2 = P_{s'} P_{s'+1} = \text{constant}$. If $\sigma (\beta_{i_1} \beta_{i_1+1}) = 0$, the last assertion still holds with a different constant. This is direct from the definitions of $a_i$ and $b_i$.  
\end{proof}

\vspace*{0.1in}
 If $\sigma (\beta_{i_k} \beta_{i_k+1}) = 1$, we set $c_k$ to be  $c_0$; otherwise we set $c_k$ to be $ c_0'$. Therefore,  the first part of lemma $\ref{002}$ says that $(d_k' \pm b_k) = c_k / (d_k \pm a_k)$, where the sign in the equation depends on the sign of $\sigma (\beta_{i_k} \beta_{i_k+1})$. \\

\begin{center}
\includegraphics[scale=0.32]{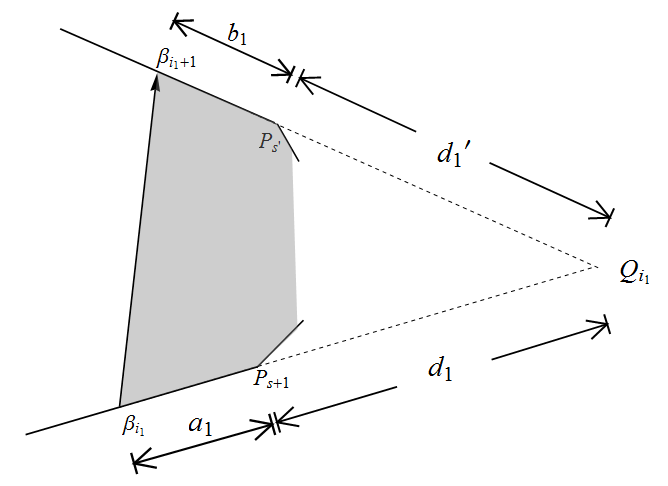} 
\includegraphics[scale=0.32]{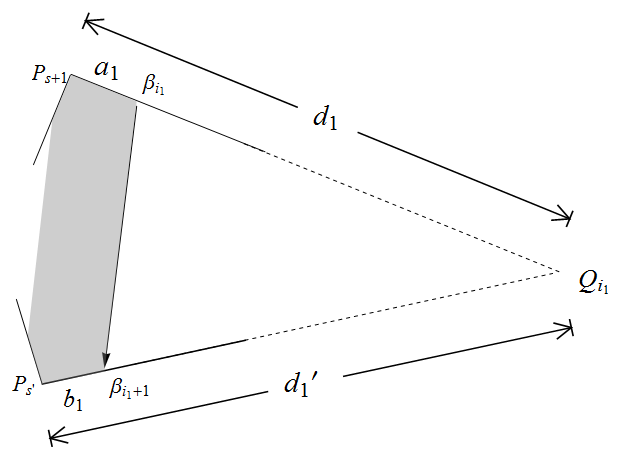} \\
\captionof{figure}{}
\label{figure06}
\end{center} 
\vspace*{0.1in}

We skip all the points that map to a parallel side in the definitions above, so we are only concerned with are those $l_i$ such that $\sigma (l_i) \ne 0$. 
We define a continuous deformation of the orbit $\mathcal{O} = \beta_0 \beta_1 ... \beta_q$ by moving $\beta_0$ with constant unit velocity.  So each $a_i(t)$ is a continuous function of $t$. Lemma $\ref{002}$ guarantees that $\dot{a_{k+1}} = - \dot{b_k}$. The reason of defining such deformation will become clear. \\

Following the discussion above we prove the following lemma:

\lemma {\label{003}} \par \indent 

If $\sigma (\beta_{i_k} \beta_{i_k+1}) = 1$, then $$\displaystyle \frac{\dot{a_{k+1}}}{\dot{a_k}} = \frac{c_k}{(a_k + d_k)^2} = \frac{b_k + d_k \,'}{a_k + d_k}.$$ 
If $\sigma (\beta_{i_k} \beta_{i_k+1}) = -1$, then  $$\displaystyle \frac{\dot{a_{k+1}}}{\dot{a_k}} = \frac{c_k}{( d_k - a_k)^2} = \frac{d_k \, ' - b_k}{ d_k - a_k}.$$

\begin{proof}
If $\sigma (\beta_{i_k} \beta_{i_k+1}) = 1$, then $\displaystyle b_k + d_k' = \frac{c_k}{a_k + d_k}$. Differentiate both sides and replace $- \dot{b_k}$ with $\dot{a_{k+1}}$ we get $\displaystyle \frac{\dot{a_{k+1}}}{\dot{a_k}} = \frac{c_k}{(a_k + d_k)^2}$. Since  $\displaystyle b_k + d_k' = \frac{c_k}{a_k + d_k}$, we obtain the desired result. 

Argument for the case when $\sigma (\beta_{i_k} \beta_{i_k+1}) = -1$ is identical. In this case $\displaystyle \frac{\dot{a_{k+1}}}{\dot{a_k}} = \frac{c_k}{( d_k - a_k)^2}$.
\end{proof}
\vspace*{0.2in}

Lemma \ref{003} also tells us that, if we move $a_1$ forward, all $a_i$ move forward, i.e., if $\dot{a_1} > 0$, then $\dot{a_i} > 0 \, \, \forall i = \{1, 2, ... q\}$. Similarly, if $\dot{a_1} < 0$, then all $\dot{a_i} < 0$. \\

We already know that if $\sigma (\beta_{i_k} \beta_{i_k+1}) = -1$, then $\sigma (\beta_{i_k+1} \beta_{i_k+2}) = 1$. Next lemma gives a stronger statement based on the results above. \\

\lemma{\label{004}} \par \indent 

Assume $\sigma (\beta_{i_k} \beta_{i_k+1}) = -1$ for some $k$. Let the orbit admit a deformation to the forward direction, so $\dot{a_i} > 0 \, \, \forall i$. Then $\displaystyle \frac{\dot{a_{k+1}}}{d_{k+1} + a_{k+1}} \geq \frac{\dot{a_k}}{d_k - a_k}$.  The equality holds if and only if $\beta_{i_k+2}$ is on the same side of $\eta$ as $\beta_{i_k}$.

\begin{proof}
We know that $\sigma (\beta_{i_k} \beta_{i_k+1}) = -1$ and $\sigma (\beta_{i_k +1} \beta_{i_k+2}) = 1$, as shown in the Figure \ref{figure07}.

\begin{center}
\includegraphics[scale=0.45]{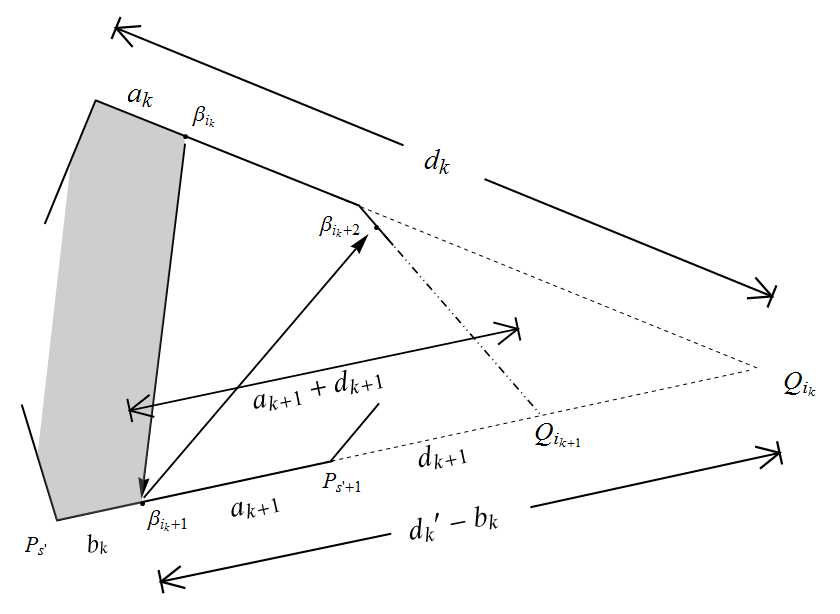}\\
\captionof{figure}{}
\label{figure07}
\end{center} 
\vspace*{0.1in}

From our definition, $P_{s'} \beta_{i_k+1} = b_k$ and $\beta_{i_k+1} Q_{i_k} = d_k \,' - b_k$. Since $\sigma (\beta_{i_k +1} \beta_{i_k+2}) \neq 0$, the side $\beta_{i_k +2} $ intersects the side $s' = P_s' P_{s'+1}$ at some point $Q_{i_{k+1}}$. Since this intersection is ``good'', i.e., $\sigma (\beta_{i_k +1} \beta_{i_k+2}) = 1$, we know that the point $Q_{i_{k+1}}$ lies on the $\mathcal{P}_{i_k 2}$ side of $\beta_{i_k +1}$ (to the right hand side of $\beta_{i_k + 1}$ in Figure \ref{figure07}). From our definitions, we know that $\beta_{i_k + 1}Q_{i_{k+1}} = a_{k+1} + d_{k+1}$. Since $\eta$ is convex, $Q_{i_{k+1}}$ has to lie between points  $Q_{i_k}$ and $\beta_{i_k +1}$, which are both on the line $P_s' Q_{i_k}$. 

It follows that $d_k\,' - b_k \geq d_{k+1} + a_{k+1}$, where the equality holds if and only if $Q_{i_k}$ coincides with $Q_{i_{k+1}}$, i.e., $\beta_{i_{k}}$ and $\beta_{i_k +2}$ are on the same sides of $\eta$. 

Since $d_k - a_k =  |\beta_{i_k} Q_{i_{k}} | > 0$, we have $$\displaystyle \frac{d_k\,' - b_k}{d_k - a_k}  \ge \frac{d_{k+1} + a_{k+1}}{d_k - a_k}.$$

Lemma \ref{003} shows that $$\displaystyle \frac{d_k\,' - b_k}{d_k - a_k} = \frac{\dot{a_{k+1}}}{\dot{a_k}}, \quad $$ so  $$\displaystyle \frac{\dot{a_{k+1}}}{\dot{a_k}} \ge \frac{d_{k+1} + a_{k+1}}{d_k - a_k}.$$

We know that, for any $ i, \, \, \, \dot{a_i} > 0$, so $$\displaystyle \frac{\dot{a_{k+1}}}{d_{k+1} + a_{k+1}} \ge \frac{\dot{a_k}}{d_k - a_k}.$$ 

From the discussion above, it is clear that the equality holds precisely when $\beta_{i_k+2}$ is on the same side of $\eta$ as $\beta_{i_k}$.
\end{proof}
\vspace*{0.2in}

Now we are ready to prove the main result, namely that $\eta$ does not contain a non-empty open interval on which all points are $q$-periodic. \\

\begin{proof}[\textbf{Proof of Theorem \ref{005}}] \indent \\

Retain our definitions of $a_1, b_1,\,  a_2,  b_2, \, ... , a_l, b_l $, the corresponding points $\beta_{i_1}, \beta_{i_2}, ..., \beta_{i_l}$, the intersections $Q_{i_1}, Q_{i_2}, ..., Q_{i_l}$ and lengths $d_1, d_1', d_2, d_2', ... ,d_l, d_l'$ from previous discussion. 

The existence of  a neighbourhood of $\beta_0$ where $\phi_q $ is identity  implies that the deformation velocity $\dot{a_{l}} = \dot{a_1}$ at $a_1(0)$. The latter statement implies $$\displaystyle \frac{\dot{a_l}}{\dot{a_{l-1}}} \frac{\dot{a_{l-1}}}{\dot{a_{l-2}}}  ... \frac{\dot{a_2}}{\dot{a_{1}}}  = 1.$$

There are two possible cases. 

(1). Assume for all $\beta_{i_k} \beta_{i_{k}+1}, \, \, k = 1, 2, ..., l, \, \, \, \, \sigma (\beta_{i_k} \beta_{i_{k+1}}) = 1$. Then by lemma \ref{003}, $$\displaystyle \frac{\dot{a_{k+1}}}{\dot{a_k}} = \frac{c_k}{(a_k + d_k)^2}.$$

So $\displaystyle \frac{\dot{a_l}}{\dot{a_{l-1}}} \frac{\dot{a_{l-1}}}{\dot{a_{l-2}}}  ... \frac{\dot{a_2}}{\dot{a_{1}}}  = 1 $ implies $$\displaystyle \frac{c_{l-1}}{(a_{l-1} + d_{l-1})^2}  \frac{c_{l-2}}{(a_{l-2} + d_{l-2})^2} ...  \frac{c_{1}}{(a_{1} + d_{1})^2}  \frac{c_{l}}{(a_{l} + d_{l})^2} = 1 $$ where $c_1, c_2 ,..., c_l $ are constants. 

This implies $$\displaystyle \prod_{k=1}^{l} (a_k + d_k)^2 = \prod_{k=1}^{l} c_k = c$$ where $c$ is a constant.  Taking logarithm we get $$\displaystyle \sum_{k=1}^{l} \ln (a_k + d_k) = \frac{1}{2} \ln c ,$$ and since $d_k$ is fixed, this implies  $$\sum_{k=1}^{l} \frac{\dot{a_k} + \dot{d_k}}{ (a_k + d_k) } = \sum_{k=1}^{l} \frac{\dot{a_k} }{ (a_k + d_k) }  = 0 .\quad $$ 

This is a contradiction since we know that all the numerators in the sum are greater than $0$ or less than 0 simultaneously, while each denominator is always greater than 0. So (1) is not possible. 

(2). There is some $k$ such that $\sigma (\beta_{i_k} \beta_{i_{k}+1}) = -1$. 

In this case,  $\displaystyle \frac{\dot{a_l}}{\dot{a_{l-1}}} \frac{\dot{a_{l-1}}}{\dot{a_{l-2}}}  ... \frac{\dot{a_2}}{\dot{a_{1}}}  = 1 $ implies $$\displaystyle    \prod_{\sigma (\beta_{i_k} \beta_{i_k+1}) = 1} \frac{c_{k}}{(a_{k} + d_{k})^2}  \prod_{\sigma (\beta_{i_k} \beta_{i_k+1}) = -1} \frac{c_k}{(d_{k} - a_k)^2} = 1 $$ where $c_1, c_2 ,..., c_l $ are constants.  This implies $$\displaystyle \prod_{\sigma (\beta_{i_k} \beta_{i_k+1}) = 1} (a_{k} + d_{k}) \, \, \, \times  \prod_{\sigma (\beta_{i_k} \beta_{i_k+1}) = -1} (d_{k} - a_k)  = c'$$ where $c'$ is a constant.

Therefore,
$$ \displaystyle \sum_{\sigma (\beta_{i_k} \beta_{i_k+1}) = 1} \ln (a_{k} + d_{k}) \, \, \, + \sum_{\sigma (\beta_{i_k} \beta_{i_k+1}) = -1} \ln (d_{k} - a_k)  = \ln c',$$ which implies
$$ \displaystyle \sum_{\sigma (\beta_{i_k} \beta_{i_k+1}) = 1} \frac{\dot{a_k}}{ (a_{k} + d_{k})} \, \, \, - \sum_{\sigma (\beta_{i_k} \beta_{i_k+1}) = -1} \frac{\dot{a_k}}{(d_{k} - a_k)  }=0.$$ 

Without loss of generality, assume $i_k = 0$, i.e., $\sigma (\beta_0 \beta_1) = -1$ where $\beta_0$ is the starting point for the $q$-periodic orbit. We know $$\sigma (\beta_1 \beta_2) = \sigma (\beta_{q-1} \beta_0)= 1$$ from lemma \ref{001} and \ref{002}.  Now consider the ordered collection of segments $$S = \{ \beta_{i_1} \beta_{i_1+1}, \beta_{i_2} \beta_{i_2+1}, ..., \beta_{i_l} \beta_{i_l+1}   \}$$  
Lemma \ref{001} and \ref{002} tell us that whenever we have $\sigma (\beta_{i_k} \beta_{i_{k}+1}) = -1$, then  $i_{k+1} = i_k + 1$, and $\sigma (\beta_{i_{k+1}} \beta_{i_{k+1}+1}) = \sigma(\beta_{i_k+1} \beta_{i_{k}+2}) = 1$. Thus we can pair up such segments $\beta_{i_k} \beta_{i_{k} +1} $ and $\beta_{i_{k}+1} \beta_{i_{k}+2}$, since the segments before and after $\beta_{i_k} \beta_{i_k+1}$ have positive signs, and the last segment in the collection, which is necessarily $\beta_{q-1} \beta_0$, also has positive sign. 
For each pair $\beta_{i_k} \beta_{i_{k} +1} $ and $\beta_{i_{k}+1} \beta_{i_{k}+2}$, we know  that $$\displaystyle \frac{\dot{a_{k+1}}}{d_{k+1} + a_{k+1}} - \frac{\dot{a_k}}{d_k - a_k} \geq 0.$$ The segments that are not paired in the collection $S$ necessarily have positive signs, so the terms related to those segments are the positive terms in the sum $$ \displaystyle \sum_{\sigma (\beta_{i_k} \beta_{i_k+1}) = 1} \frac{\dot{a_k}}{ (a_{k} + d_{k})} \, \, \, - \sum_{\sigma (\beta_{i_k} \beta_{i_k+1}) = -1} \frac{\dot{a_k}}{(d_{k} - a_k)  }=0.$$

If not all segments in $S$ are paired up as above, then $$ \displaystyle \sum_{\sigma (\beta_{i_k} \beta_{i_k+1}) = 1} \frac{\dot{a_k}}{ (a_{k} + d_{k})} \, \, \, - \sum_{\sigma (\beta_{i_k} \beta_{i_k+1}) = -1} \frac{\dot{a_k}}{(d_{k} - a_k)  } > 0,$$ and this leads to contradiction. 

Now assume all segments are paired up, so $l$ is even and the signs of the segments in $S$ alternate as $\{( -, +) , ( -, +) ,  ... , ( -, +) \}$. Again, since $\displaystyle \frac{\dot{a_{k+1}}}{d_{k+1} + a_{k+1}} - \frac{\dot{a_k}}{d_k - a_k} \geq 0$ for each pair, we have  $$ \displaystyle \sum_{\sigma (\beta_{i_k} \beta_{i_k+1}) = 1} \frac{\dot{a_k}}{ (a_{k} + d_{k})} \, \, \, - \sum_{\sigma (\beta_{i_k} \beta_{i_k+1}) = -1} \frac{\dot{a_k}}{(d_{k} - a_k)  } \geq 0,$$  and the equality holds if and only if all $\beta_{i_k}$ and $\beta_{i_k +2} = \beta_{i_{k+2}}$ are on the same side of $\eta$. This leads to the conclusion that $\beta_{i_1}, \beta_{i_3}, \beta_{i_5}, ... , \beta_{i_{l-1}}$ are all on the same side as $\beta_{i_1}$. Furthermore, since the area we cut off from $\eta$ is strictly less than a half of the total area, the points $\beta_{i_1}, \beta_{i_3}, \beta_{i_5}, ... , \beta_{i_{l-1}}$ are ordered as listed on the line with no points coinciding with the other. This suggests that $\beta_{i_{l+1}}$  does not coincide with $\beta_{i_1}$, which is clearly a contradiction. So the theorem is proved.
\end{proof}

\indent \\
In the proof of theorem \ref{005} we assumed that $\mathcal{O}$ does not contain any corner of $\eta$. Now we finish the proof of the main theorem by showing that the result holds when $\mathcal{O}$ does contain corners. We can argue by contradiction. Assume that $\mathcal{O} $ touches some corners $p_{k_1}, ..., p_{k_m}$ of $\eta$ and the $q^{th}$ iteration of the billiard map is locally identity.  Without loss of generality, let $\beta_0 = p_{k_1}$, so there is an open neighbourhood $J_0$ of $\beta_0$ such that the $q^{th}$ iteration map is identity. Since the orbit has finite period $q$, there exists some $\beta_0' \in J_0$ such that $\beta_0'$ also leads a $q$-periodic orbit $\mathcal{O}'$ (since $\beta_0' \in J_0$) and the new orbit  $\mathcal{O}'$ does not contain corners. (We just need to perturb the orbit slightly). This contradicts theorem $\ref{005} $ that we just proved. This finishes the proof of our main theorem stated in the beginning of the article. \\

Therefore, by the discussion in part 1 and the main theorem , we obtain:

\theorem{\label{007}} The rotation number of the circular homeomorphism induced from the area construction of any convex polygon is always a devil's staircase function of the area parameter. \\

We then have an interesting corollary, compared to results by Innami \cite{IN} and Baryshnikov / Zharnitsky \cite{BY}.

\corollary{\label{008}} A circular map induced from area construction of any convex polygon cannot consist solely of periodic points. That is to say, considering a convex polygonal invariant curve and a piecewise hyperbolic table resulted from the area construction, the convex polygonal invariant curve contains non-periodic points under the corresponding outer billiard map. \\

\remark{\label{009}}  Finally we remark that, for our result to hold, the area used in the area construction need not to be fixed. We could consider a generalized map defined on sides of polygons (locally). We still use the area construction, but instead of cutting off area $A$, we cut off areas $A_1, A_2 ,..., A_q$ each time in our construction to define a sequence of $q+1$ points $\gamma_0, ..., \gamma_{q}$ where $\gamma_0 = \gamma_{q}$. We call this orbit a fake periodic orbit. Then locally we can define ``billiard  map" on a small interval around $\gamma_0$, by the area construction of cutting area $A_i$ to obtain $\gamma_i$. Then the methods we use to prove the main theorem still apply and we conclude that the $q^{th}$ iteration of the``billiard map " in this case cannot be identity on the interval. \\ \\

\newpage

Now we present several numerically generated graphs of the devil's staircase functions. Figure \ref{figure08} shows the simplest case, where the polygon is a square.

\begin{center}
\includegraphics[scale=0.6]{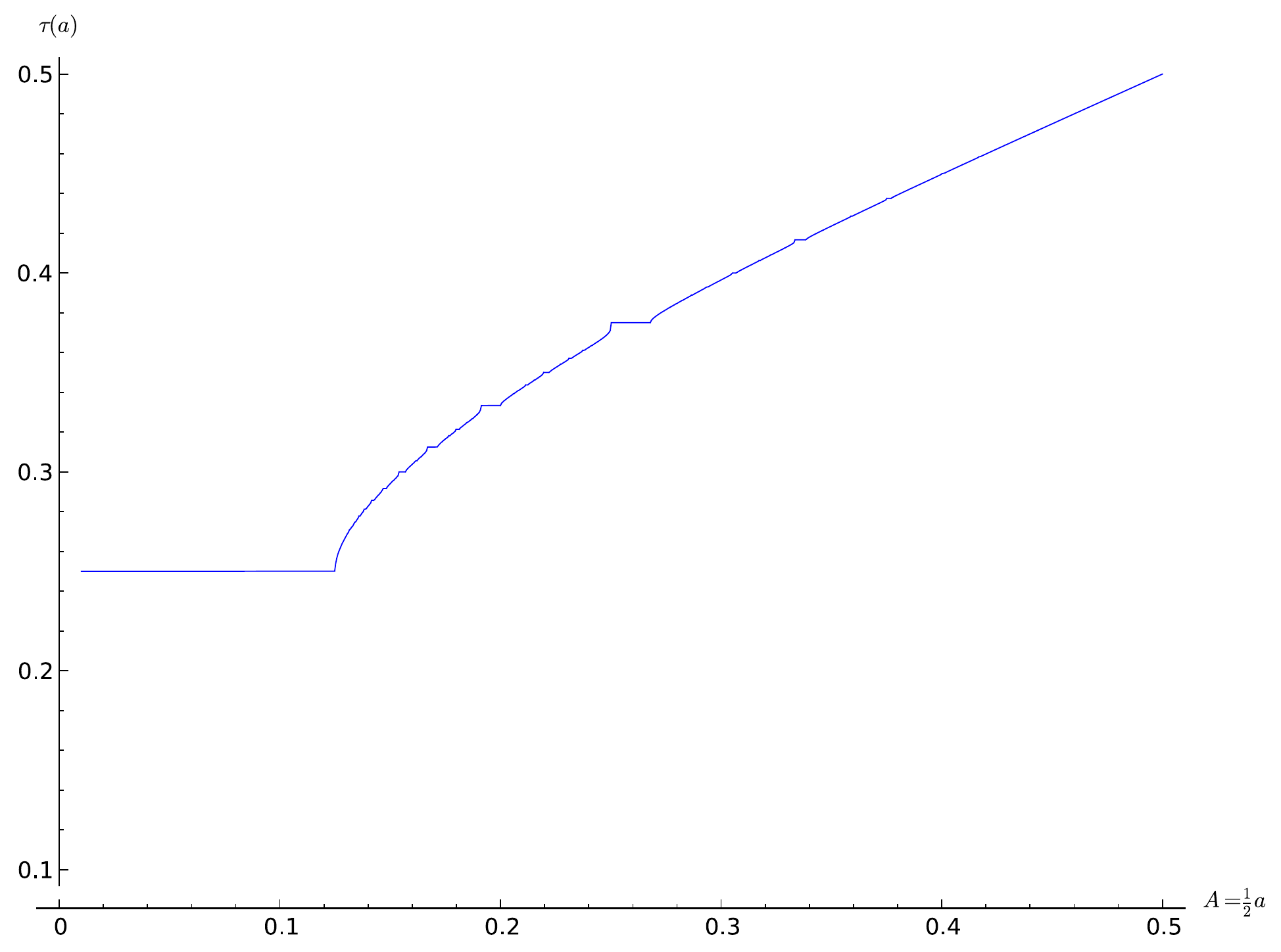}\\
\captionof{figure}{ $\eta$  is a square}
\label{figure08}
\end{center} 
\vspace*{0.3in}

A small part of Figure \ref{figure08} is zoomed in to show the detailed features of the devil's staircase behaviour.  The domain of the $x$-axis of Figure \ref{figure11} is  $ (0.13, 0.21)$.

\begin{center}
\includegraphics[scale=0.5]{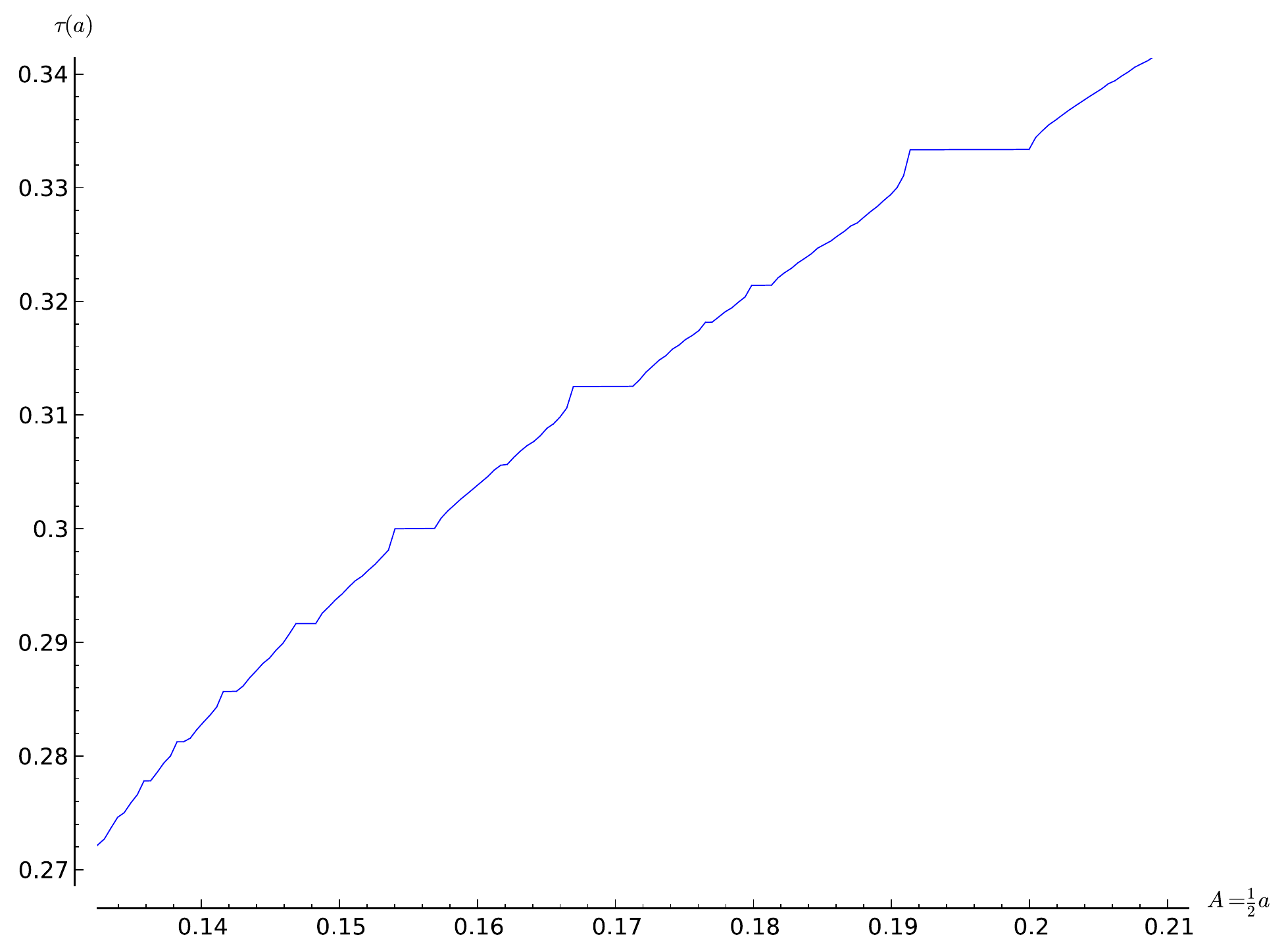}\\
\captionof{figure}{ A zoomed-in portion of Figure \ref{figure08}}
\label{figure11}
\end{center} 

Here we present the devil's staircases for a regular pentagon and for an irregular pentagon with coordinates $\{(0,0), (2,0), (2,1), (1,2), (0,1) \}$: 

\begin{center}
\includegraphics[scale=0.6]{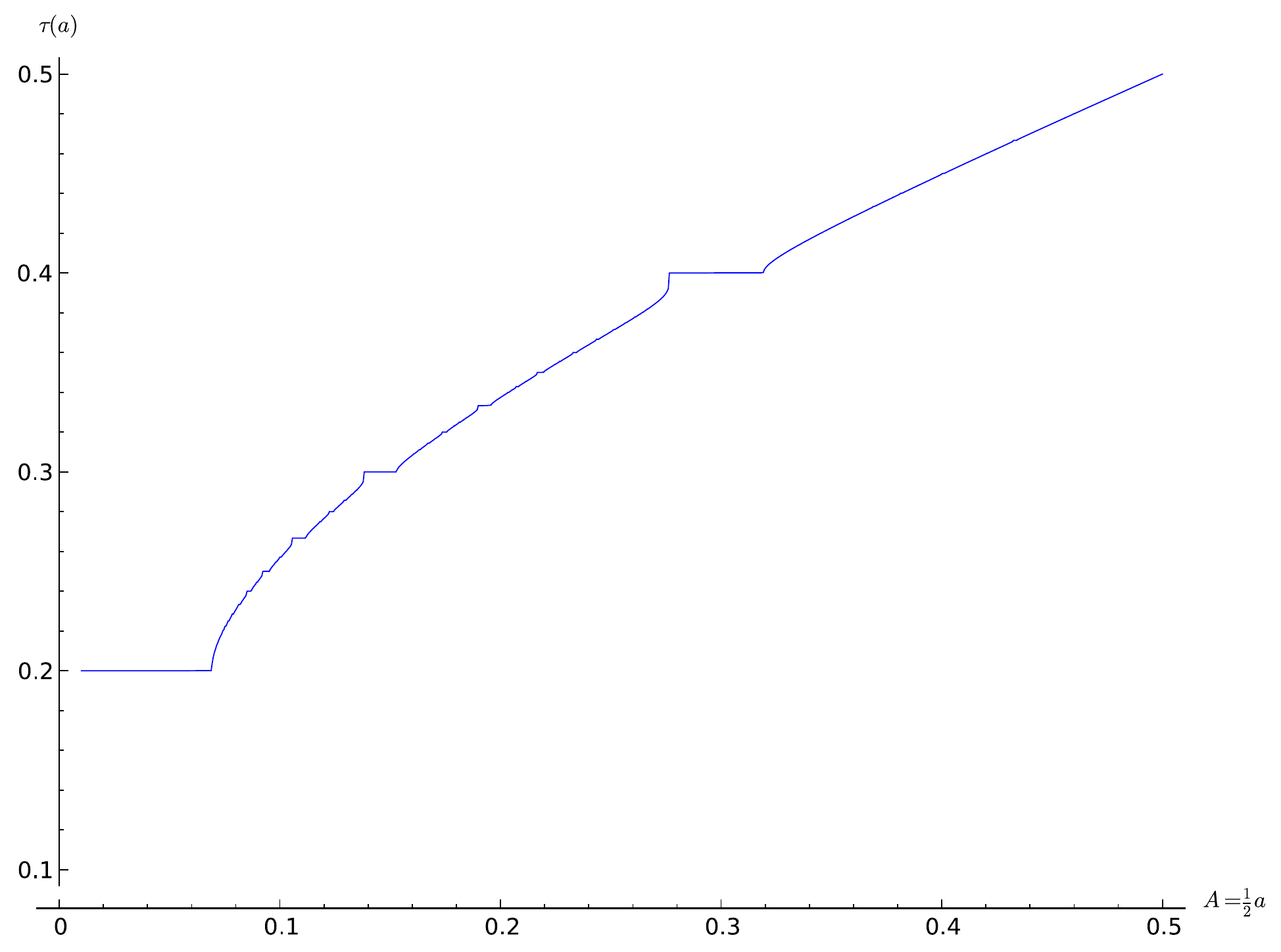}\\
\captionof{figure}{ $\eta$  is a regular pentagon}
\label{figure09}
\end{center} 
\vspace*{0.3in}

\begin{center}
\includegraphics[scale=0.6]{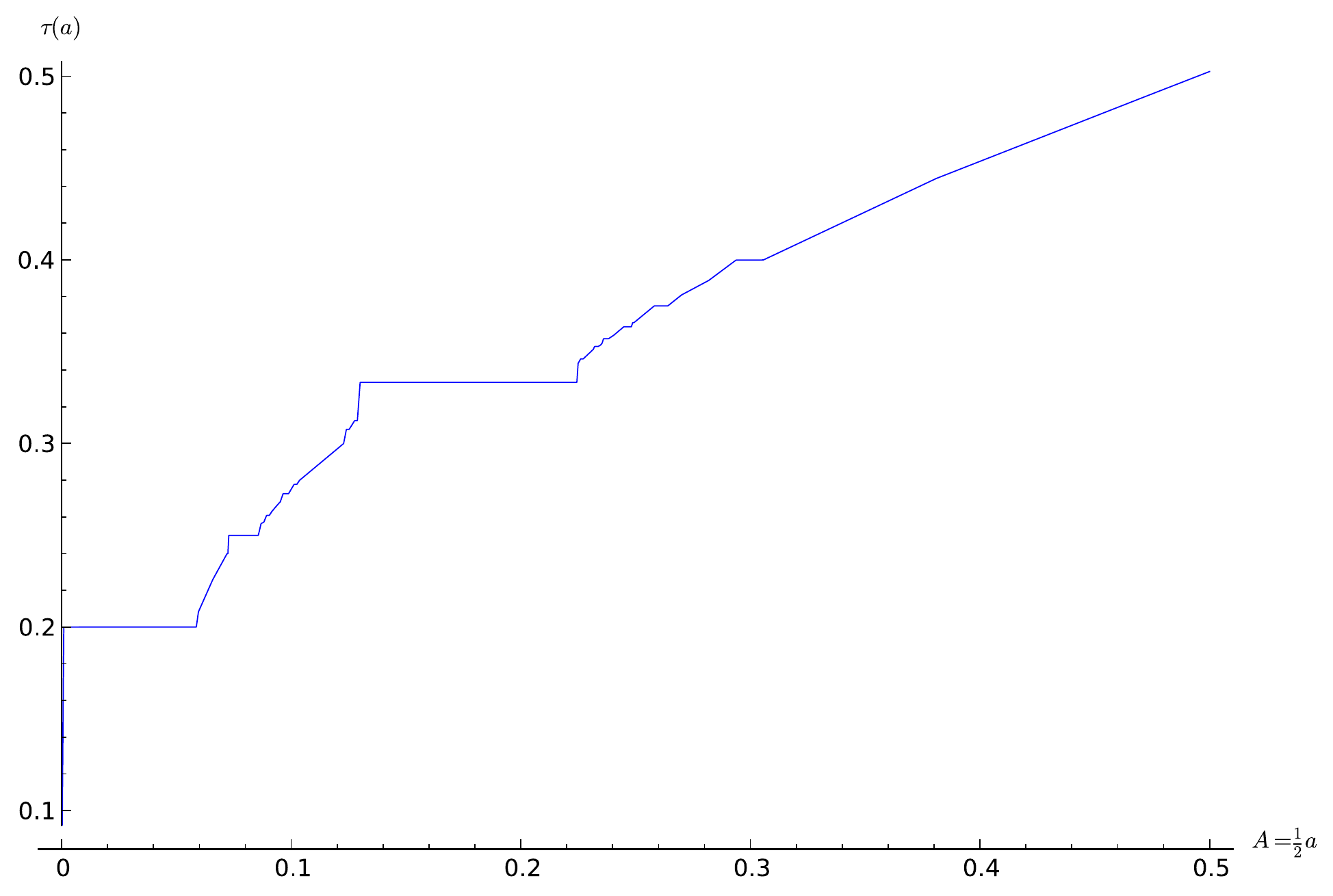}\\
\captionof{figure}{ $\eta$ is a pentagon with coordinates $\{(0,0), (2,0), (2,1), (1,2), (0,1) \}$}
\label{figure10}
\end{center}

\newpage

\section{ \fontsize{15}{14}\selectfont \textit{Discussion on required hypothesis of the result and open questions}}

Consider the following setup: given an ordered collection of lines $l_1, l_2, ..., l_n$, not necessarily distinct, in Euclidian plane $\mathbb{R}^2$, define a series of functions $f_k: l_k \rightarrow l_{k+1}, \, \, \, \, k = 1, ..., n ,\, \, \,$ with $n+1$ set to be $1$.  Each function is given by the associated area construction map between two consecutive lines $l_k, l_{k+1}$.  This map is essentially a map from $\mathbb{R}P^{1} \rightarrow \mathbb{R}P^{1}$, sending $\infty \in l_k$ to the intersection $l_k \cap l_{k+1} = I_k$ of the two lines, and $I_k \in l_k$ to $\infty \in l_{k+1}$. \\

We ask the following questions: does there exist such a collection of ordered lines with the defined maps $f_1, ..., f_n$ such that $f_n \circ f_{n-1} \circ ... \circ f_1 = identity$? Is the convex property of the polygon $\eta$ essential for theorem \ref{005} to hold? We ask for the most broad generalization possible.\\

The answer to the first question is yes, here we present two simple examples that satisfy the requirements.

\begin{center}
\includegraphics[scale=0.3]{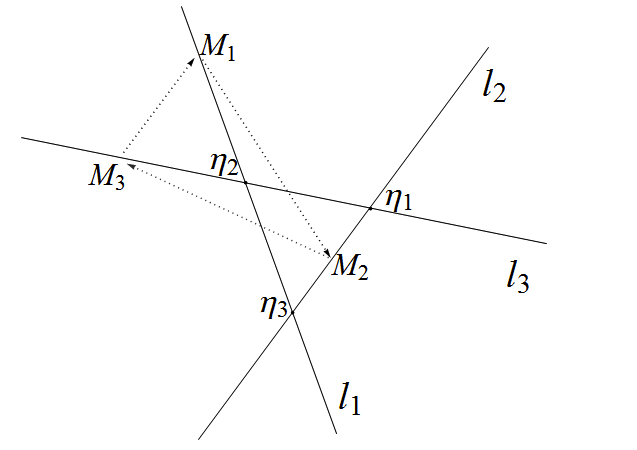}
\includegraphics[scale=0.3]{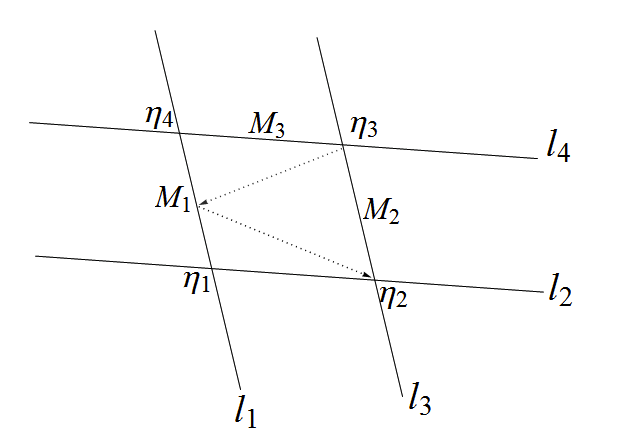}  \\
\captionof{figure}{}
\label{figure12}
\end{center}  
\vspace*{0.2in}

First consider a collection of three ordered lines $l_1, l_2, l_3$ that are not concurrent, as shown in Figure \ref{figure12}. $\eta_1, \eta_2, \eta_3$ are the three intersections of the three lines respectively, $M_1 \in l_1, M_2 \in l_2$, and $ M_3 \in l_3 $. Let $M_1 \eta_2 = \eta_2 \eta_3$, $M_3 \eta_2 = \eta_2 \eta_1$ while $M_2$ is the midpoint of $\eta_1 \eta_2$. We define three maps $f_1, f_2, f_3$ from area construction on consecutive lines by cutting the whole area of the enclosed triangle $\eta_1 \eta_2 \eta_3$. \\

\proposition{\label{10}} For the collection of $l_1, l_2, l_3$ and maps $f_1, f_2, f_3$ defined above, $f_3 \circ f_2 \circ f_1 : l_1 \rightarrow l_1$ is identity. 
\begin{proof}

We give each of the three lines Euclidean coordinates. It does not matter where we set the origins to be. We claim each map $f_i$ is a M\"{o}bius transformation on the coordinates of the lines. The area construction guarantees that $M_i \eta_{i+2} \cdot M_{i+1} \eta_{i+2} =  \text{constant} $.  Thus $M_{i+1} \eta_{i+2} = {\text{constant}}/{M_{i} \eta_{i+2}}$, which is indeed a M\"{o}bius transformation when we use the defined coordinates.  The composition $f_3 \circ f_2 \circ f_1$ is therefore also a M\"{o}bius  transformation, so we only need to prove that it fixes three distinct points to show it is identity.  \\

First consider $\infty \in l_1$. We have $f_1(\infty) = \eta_3, \, \, f_2(\eta_3) = \eta_2$ by cutting the area of the triangle, and $f_3(\eta_2) = \infty \in l_1$.  \\

Now consider point $\eta_3 \in l_1$. We have $f_1(\eta_3 ) = \infty \in l_2, \, \, \, f_2(\infty) = \eta_1 \in l_3, \,\, \, \text{and } f_3(\eta_1) = \eta_3$.  \\ 

Finally, consider $M_1 \in l_1$. Clearly the area of $\triangle M_1 M_2 \eta_3$ is the same as the area of $\triangle \eta_1 \eta_2 \eta_3$, which is the same as the area of $\triangle M_3 M_2 \eta_1$, so $f_1(M_1) = M_2, \, \, f_2(M_2) = M_3, \, \, f_3(M_3) = M_1$. \\

Therefore, $\infty, M_1, \,\, \eta_3$ are three fixed point of the M\"{o}bius  transformation $f_3 \circ f_2 \circ f_1$, and the composed map is identity. 
\end{proof}

Now consider a collection of four ordered lines $l_1, l_2, l_3, l_4$, where $l_1 \parallel l_3, \, \, l_2 \parallel l_4 $. The enclosed area forms a parallelogram $\eta_1 \eta_2 \eta_3 \eta_4$, whose vertices are intersections of pairs of consecutive lines. Let $M_1, M_2, M_3$ be the midpoints of three of sides of the parallelogram as shown in Figure \ref{figure12}. Define $g_1, g_2, g_3, g_4$ analogously by cutting a quarter of the area of the parallelogram. 

\proposition{\label{11}}For the collection of $l_1, l_2, l_3, l_4$ and maps $g_1, g_2, g_3, g_4$ defined above, $g_4 \circ g_3 \circ g_2 \circ g_1 : l_1 \rightarrow l_1$ is identity. 
\begin{proof}
Similar to the proof of proposition \ref{10}, we keep track of the points $0, \infty $ and midpoints $M_1$. \\

$M_1 \in l_1  \, \, \, \xrightarrow{ \, \, \, \, f_1 \, \, \, \,  } \, \,  \eta_2 \in l_2  \, \, \, \xrightarrow{ \, \, \, \, f_2 \, \, \, \,  }  \, \, \infty \in l_3  \, \, \, \xrightarrow{ \, \, \, \, f_3 \, \, \, \,  } \, \, \eta_3 \in l_4 \,  \, \, \xrightarrow{ \, \, \, \, f_4 \, \, \, \,  }  \, \, M_1 \in l_1   $

$\infty \in l_1  \, \, \, \xrightarrow{ \, \, \, \, f_1 \, \, \, \,  } \, \,  \eta_1 \in l_2  \, \, \, \xrightarrow{ \, \, \, \, f_2 \, \, \, \,  }  \, \, M_2 \in l_3  \, \, \, \xrightarrow{ \, \, \, \, f_3 \, \, \, \,  } \, \, \eta_4 \in l_4 \,  \, \, \xrightarrow{ \, \, \, \, f_4 \, \, \, \,  }  \, \, \infty \in l_1   $

$\eta_1 \in l_1  \, \, \, \xrightarrow{ \, \, \, \, f_1 \, \, \, \,  } \, \,  \infty \in l_2  \, \, \, \xrightarrow{ \, \, \, \, f_2 \, \, \, \,  }  \, \, \eta_2 \in l_3  \, \, \, \xrightarrow{ \, \, \, \, f_3 \, \, \, \,  } \, \, M_3 \in l_4 \,  \, \, \xrightarrow{ \, \, \, \, f_4 \, \, \, \,  }  \, \, \eta_1 \in l_1   $

So the composed map is identity.
\end{proof}
\vspace*{0.15in}

This tells us that our theorem cannot be generalized to arbitrary collection of lines. \\

\textnormal{\textbf{Open Questions}}.\\

 The following questions are still open:\\

1. Does there exist non-convex simple polygons such that some $n^{th}$ iteration map gives identity?  That is, for simple polygons, is convexity a required condition?  \\

2. Given a convex smooth ($C^1$ or $C^\infty$) simple closed curve, we can still construct the corresponding table such that the curve is an invariant curve of the outer billiard system. We want to know to what extent the result of Theorem \ref{005} still holds. \\

We know that if the curve is a circle, then the corresponding table is also a circle. In which case, the outer billiard map is simply a rotation map, while the rotation number is a smooth strictly increasing function of the area parameter. Note that the problem for an ellipse is the same as for the circle since the problem is affine-invariant. It is also proved in \cite{GD} by Genin and Tabachnikov that there exist non-circular outer billiards having invariant curves consisting of periodic points. They constructed such curves from perturbing a circle. Therefore, for these two families of curves, our main result does not hold. The following question is still open:  for which curves can we conclude that the rotation number is a devil's staircase of the area parameter?

\newpage
\section{ \fontsize{15}{14}\selectfont \textit{Acknowledgements}}

This result was obtained during the summer@ICERM program 2013, and was motivated by questions asked by  Sergei Tabachnikov. I would like to thank the advisors Sergei Tabachnikov, Ryan Greene and Diana Davis for valuable discussions. In addition, Ryan Greene provided the codes in Sage to numerically generate graphs of the devil's staircase functions. The author is funded by UTRA awards from Brown University, under supervision of  Thomas Banchoff. Finally, I want to thank ICERM for its consistent support and the opportunity to carry out this work. \\ \\ \\

\end{document}